\documentclass[12pt,reqno,a4paper]{amsart}
\usepackage{
    amsmath,  amsfonts, amssymb,  amsthm,   amscd,
    gensymb,  graphicx, comment,  etoolbox, url,
    booktabs, stackrel, mathtools,enumitem, mathdots,  microtype, lmodern,    mathrsfs, graphicx, tikz,  longtable,tabularx, float, tikz, pst-node, tikz-cd, multirow, tabularx, amscd,  bm, array, makecell, diagbox, booktabs,ragged2e, caption, subcaption }
\usepackage{makecell,slashbox}
\usepackage{algorithm}
\usepackage{algpseudocode}
\usepackage{esvect}

\MakeRobust{\vv}
\usepackage{xcolor}
\usepackage[utf8]{inputenc}
\usepackage{microtype, fullpage, wrapfig,textcomp,mathrsfs,csquotes,fbb}
\usepackage[colorlinks=true, linkcolor=blue, citecolor=blue, urlcolor=blue, breaklinks=true]{hyperref}
\usepackage[capitalise]{cleveref}
\usepackage{todonotes}
\usetikzlibrary{positioning}
\usetikzlibrary{shapes,arrows.meta,calc}

\usetikzlibrary{arrows}

\newtheorem{theorem}{Theorem}[section]

\newtheorem{corollary}[theorem]{Corollary}
\newtheorem{definition}[theorem]{Definition}
\newtheorem{example}[theorem]{Example}
\newtheorem{lemma}[theorem]{Lemma}

\newtheorem{proposition}[theorem]{Proposition}
\newtheorem{remark}[theorem]{Remark}

\newcommand\R{\mathbb{R}}

\newcommand{\tc}{\vv{\mathrm{TC}}}
\newcommand{\cat}{\vv{\mathrm{cat}}}

\newcommand{\TC}{\mathrm{TC}}

\newcolumntype{x}[1]{>{\centering\arraybackslash}p{#1}}

\begin{document}
\title{On the topological complexity of directed parametrized motion planning}

\author[S. Datta ]{Sutirtha Datta}
\address{Department of Mathematics, Indian Institute of Science Education and Research Pune, India}
\email{sutirtha2702@mail.com}
\author[N. Daundkar ]{Navnath Daundkar}
\address{Department of Mathematics, Indian Institute of Technology Madras, Chennai, India}
\email{navnath@iitm.ac.in}
\author[A. Sarkar ]{Abhishek Sarkar}
\address{Department of Mathematics, Indian Institute of Science Education and Research Pune, India}
\email{abhisheksarkar49@gmail.com}


\begin{abstract}
We introduce and study a parametrized analogue of the directed topological complexity, originally developed by Goubault, Farber, and Sagnier. We establish the fibrewise basic dihomotopy invariance of directed parametrized topological complexity and explore its relationship with the parametrized topological complexity. In addition, we introduce the concept of the directed Lusternik–Schnirelmann (LS) category, prove its basic dihomotopy invariance, and investigate its connections with both directed topological complexity and directed parametrized topological complexity. We further investigate additional properties of our invariant and examine its connections with several other invariants that arise naturally in the context of topological robotics. Moreover, we compute the directed parametrized topological complexity of the Hopf fibrations and the Fadell-Neuwirth fibrations having specific directed fibration structures.
\end{abstract}
\maketitle
\maketitle

\section{Introduction}\label{sec:intro}
A  \emph{motion planning algorithm} on a path-connected topological space $X$ is defined  as a section of the free path space fibration $\pi\colon X^I\to X\times X,$ defined by
$\pi(\gamma)=(\gamma(0),\gamma(1)), $
where  $X^I$ denotes the free path space of $X$, equipped with the compact open topology.
To analyze the complexity of designing a motion planning algorithm for the configuration space $X$ of a mechanical system, Farber \cite{FarberTC} introduced the notion of topological complexity. 
The \emph{topological complexity} of a space $X$, denoted by $\TC(X)$, is defined as the smallest natural number $r$ for which $X\times X$ can be covered by open sets $\{U_1,\dots, U_r\}$, with each $U_i$ admitting a continuous local section of $\pi$. 
The number $\TC(X)$ represents the minimal number of continuous rules required to implement a motion planning algorithm in the space $X$.
Farber \cite[Theorem 3]{FarberTC} showed that $\TC(X)$ is a numerical homotopy invariant of a space $X$. 

 The notion of directed topological complexity for Euclidean neighborhood retracts (ENR) was introduced more recently (see \cite{goubault2017directed}, \cite{goubault2020directed}), which provides a natural framework for studying motion planning problems where control constraints restrict the movements of different components of a mechanical system.
 The origin of the directed algebraic topology dates back to the foundational works of Grandis \cite{grandis2002directed, grandis2003directed, grandis2009directed}. Directed topology has found a wide range of applications in concurrency theory, hybrid dynamical systems, and related fields (see \cite{fajstrup2016directed}).
In \cite{goubault2017directed}, Goubault defined a variant of topological complexity for \textit{directed topological spaces} and showed that it is invariant under a suitable notion of directed homotopy equivalence. Later, in \cite{goubault2020directed}, he extended the theory with Farber and Sagnier. Subsequently, Borat and Grant studied the directed topological complexity of directed spheres with the directed structure induced by the directed cube on their boundaries in \cite{borat2020directed}. It is worth noting that computations of directed topological complexity are extremely challenging, and known results are available only for specific cases such as directed spheres and directed graphs.

A \textit{directed space} (or \textit{d-space}) is a space $X$ equipped with a distinguished class of paths $dX$ in $X$, called \textit{directed paths} (or \emph{d-paths}), which satisfy certain axioms (see Definition \ref{d-space}). The directed paths of a d-space $X$ form a subspace  
$dX  \subseteq X^I.$ 
The \textit{free path space fibration} restricts to a map  
\[
\vv{\pi}\colon dX \to \Gamma_X,
\]  
where \( \Gamma_X \subseteq X \times X \), is the set of pairs \( (x, y) \in X \times X \) such that there exists a directed path from \( x \) to \( y \). A \textit{directed motion planner} on a subset \( A \subseteq \Gamma_X \) is a continuous local section of \( \vv{\pi} \) on \( A \). The \textit{directed topological complexity} of the d-space \( X \), denoted by $\tc(X)$, is the smallest natural number \( k \) such that $\Gamma_X$ can be partitioned into \( k \) disjoint Euclidean Neighborhood Retracts (ENRs), each admitting a directed motion planner.

In \cite{goubault2020directed}, the authors introduced an appropriate notion of homotopy equivalence, called basic dihomotopy equivalence, in directed algebraic topology and showed that the directed topological complexity is invariant with respect to this equivalence. Furthermore, they established several properties of $\tc(X)$, including a product formula, and provided examples where the directed topological complexity can be estimated.
For strongly connected directed spaces, they proved in \cite[Proposition 2]{goubault2020directed} that $\TC(X) \leq \tc(X)$. In general, however, there is no relation between $\tc$ and $\TC$, as the authors provided examples of directed spaces for which these two invariants differ. Thus, $\tc$ stands out as a completely new invariant, which may be independently interesting to study and explore.

\subsection{Directed LS category} 
There is an old invariant called LS category, a close relative of topological complexity, which was introduced by Lusternik and Schnirelmann in \cite{LScat}. The \emph{LS category} of a space $X$ is denoted by $\mathrm{cat}(X)$, which is the least natural number of open subsets that cover $X$ such that the inclusion on each open set is nullhomotopic.
In \cite{FarberTC}, Farber proved the following famous inequality
$\mathrm{cat}(X)\leq \TC(X)\leq \mathrm{cat}(X \times X).$
To introduce an analogous notion for a d-space \( X \), called the \emph{directed LS category} of \( X \), denoted by \(\cat(X)\), we proceed as follows. Consider a directed space $X$ with an initial point
  $x_0$, i.e., there is a directed path from $x_0$ to every point in $X$. For a directed space $X$ with an initial point $x_0$, let $d_0X$ be the set of directed paths starting at the initial point $x_0$. The evaluation map $e_X\colon d_0X\to X$ is defined by sending a d-path to its endpoint.
 Unlike the undirected setting, for a path-connected \(d\)-space \(X\), the existence of a continuous section of $e_X\colon d_0X \to X,$
is \emph{not} equivalent to the dicontractibility of \(X\). We give a brief indication of the strategy for defining the directed LS-category in Remark \ref{motivation for directed cat}.
Here, $\cat(X)$ is defined to be the smallest natural 
number $k$ such that $X$ can be covered by $k$ ENRs $\{U_i\} _{i=1}^{k}$, where each $U_i$ admits a continuous directed section of the map $e_X$.
We show the following identities for a d-space $X$ with initial point:
    $$\mathrm{cat}(X) \leq \cat(X) \text{~and~}~\cat(X)\leq \tc(X).$$ 
For a d-space $X$ with an initial point, we establish the basic dihomotopy invariance of $\cat(X)$. In Example \ref{d-sphere}, we compute $\cat(S^n)$ for the $n$-dimensional sphere $S^n$ with a standard directed structure and with the origin $\bar{0} \in \mathbb{R}^{n+1}$ as initial point, for $n \in \mathbb{N}$.
\subsection{Directed parametrized topological complexity}
The notion of parametrized topological complexity, introduced by Cohen, Farber, and Weinberger \cite{C-F-W}, extends motion planning to settings with external parameters. A parametrized motion planning algorithm assigns, to each pair of configurations under the same external condition, a continuous path that respects and remains constant with respect to that condition. 
We now recall the notion of parametrized topological complexity in brief. For a fibration $p\colon  E \to B$, 
the space $E^I_B$ denotes the space of all paths in $E$ with image in a single fibre. The restriction of the free path space fibration on $E^{{I}}_B \subseteq E^I$ produces a fibration ${\Pi}\colon E^{{I}}_B\to {E \times_B E}$.
The \emph{parametrized topological complexity} of a fibration $p\colon E \to B$ denoted by $\TC[p\colon E\to B]$  is the smallest natural number $k$ such that there is an open cover $\{U_1,\dots, U_k\}$ of $E\times_B E$, where each open set $U_i$ admits a continuous section of $\Pi$. 
The reader is referred to \cite{C-F-W}, \cite{PTCcolfree}, \cite{ptcspherebundles} for several interesting results related to parametrized topological complexity. Moreover, the notion of parametrized topological complexity of fibrations is extended to fibrewise spaces in \cite{GC} by Garc\'{\i}a-Calcines.

 In this paper, we unify the notions of directed topological complexity and parametrized topological complexity by introducing the concept of directed parametrized topological complexity. In applications, this invariant measures the minimal number of continuous motion-planning rules that simultaneously respect directionality (e.g., non-reversibility or causal structure) and external parameters (such as obstacle configurations). It thereby provides a refined tool for evaluating and comparing the intrinsic complexity of control-constrained motion-planning tasks.

Suppose $E$ and $B$ are ENRs. Then,
 for a d-fibration $p\colon  E \to B$, we
 denote the space of d-paths in $E$ with image in a single fibre by $dE_B$. Consider the  restriction of $\Pi$ to the subspace $dE_B$ of $E^I_B$ with image as the space of end points $\Gamma_{E, B}$ of the d-paths in $dE_B$. We denote this restriction by $${\vv{\Pi}}\colon dE_B\to {\Gamma_{E, B}}.$$
The \emph{directed parametrized topological complexity} of a d-fibration $p \colon E \to B$, denoted by $\vv{\TC}[p\colon E\to B]$, is the smallest natural number $k$ such that there is a cover by ENRs $\{U_1,\dots, U_k\}$ of $\Gamma_{E, B} \subseteq E\times_B E$, where each ENR $U_i$ admits a continuous section of $\vv{\Pi}$. 
If the fibre $F$ of a d-fibration $p\colon E\to B$ is strongly connected, then we show that
$$ \mathrm{TC}[p\colon E\to B]\leq \tc[p\colon E\to B].$$ 
Moreover, we show that these invariants may differ by providing Example~\ref{undirected TC < directed TC}, considering a disc bundle equipped with a directed structure on the $2$-dimensional disc $D^2$, for which $\tc(D^2)=2$, in contrast to the usual $\mathrm{TC}(D^2)=1$. In general, no comparison holds between them; to demonstrate this, we construct Example~\ref{directed TC< undirected TC}, considering a circle bundle equipped with a directed structure on the unit circle $S^1$, for which $\tc(S^1)=1$, in contrast to the usual $\mathrm{TC}(S^1)=2$.

We introduce the fibrewise version of basic dihomotopy equivalence for d-fibrations and show that the directed parametrized topological complexity is a numerical fibrewise basic dihomotopy invariant. 
For a fibration $p\colon E \rightarrow B$ with a d-structure on $B$, we associate a natural d-structure on $E$ for which the fibration $p$ becomes a directed fibration. This directed structure allows us to establish the following equality (\cref{directed-usual equality})
\[\tc[p\colon E\to B] = \mathrm{TC}[p\colon E\to B].\]
 The problem of obstacle-avoiding, collision-free motion of multiple robots in the presence of multiple obstacles with unknown positions was studied in \cite{C-F-W, PTCcolfree}. The key ingredients were the \emph{Fadell-Neuwirth fibrations} and their parametrized topological complexity. In this work, we compute a \emph{directed} version of this invariant using the results described above.

\subsection{Layout of the article}
In \cref{sec : background}, we recall key notions from directed algebraic topology that are essential for the subsequent sections. We begin with the definitions of directed spaces, directed maps, and directed homotopies. We also review the concept of basic dihomotopy equivalence and its relevance to directed topological complexity. Additionally, we recall certain special spaces, such as strongly connected spaces and regular d-spaces.

In \cref{sec: directed LS and directed TC}, we introduce the notion of directed LS category and establish it as a numerical invariant under the basic dihomotopy equivalence relation in \cref{thm: dihtopy equivalent then equal cat}. Next in Proposition \ref{prop: cat ineq & cat 1} and Proposition \ref{prop: rel cat-tc} we relate the newly defined directed LS category with the notions of LS category, directed topological complexity and topological complexity.  

Finally in \cref{Sec 4: directed parametrized topological complexity}, we define directed parametrized topological complexity as a fusion of the notions of parametrized topological complexity and directed topological complexity. We divided this section into two subsections, \Cref{properties} and \Cref{invariance}, one is to study its properties and another to study its invariance. 

In \cref{properties}, we establish various relationships among existing notions of topological complexities with this new notion. 
In Proposition \ref{prop:ineq: pullback and tc}, for a d-fibration $p\colon  E \to B$, we establish an inequality relating $\tc$ of the fibre and $\tc[p\colon  E \to B]$. In particular, for a trivial d-fibration, Proposition \ref{prop: trivial d-fib ; equality of tc} provides an equality. In \cref{thm: dicontractibility implies tc 1}, under mild assumption on the fibre, we show that for a d-fibration $p\colon  E \to B$, having dicontractible fibre is equivalent to obtaining $\tc[p\colon  E \to B]=1$. In Proposition \ref{prop: dir and undir paramtrized TC inequality}, we establish an inequality relating the directed and undirected versions of parametrized topological complexity. Next, in Proposition \ref{prop: product ineq dir parametrized tc}, under the regularity assumption on the d-fibrations, we establish a product inequality for directed parametrized topological complexity. In Proposition \ref{prop: paramtrized tc leq directed LS cat}, for a d-fibration $p\colon  E \to B$ we establish an inequality among the parametrized topological complexity of $p$ and the directed LS category of $E \times_B  E$.

 In \cref{invariance}, we introduce the notion of fibrewise basic dihomotopy equivalence in Definition \ref{def: d-fib hteq} and  establish the equality for $\tc[p\colon  E \to B]$ and $\tc[p'\colon E' \to B]$ for two fibrewise basic dihomotopy equivalent d-fibrations $p$ and $p'$ in \cref{thm: basic dihomotopy eqv fib have same TC}. 

In \cref{Sec 5: computations}, we consider a particular kind of non-trivial but very natural directed structure on a fibration and using that we  compute the directed parametrized topological complexity for the Hopf fibrations and the Fadell-Neuwirth fibrations.

\section{Background for directed algebraic topology}\label{sec : background}
We first recall the definitions of directed spaces, directed paths, and directed maps, which we refer to as d-spaces, d-paths, and d-maps, respectively. We then review the notions of basic dihomotopy equivalence and regular d-spaces, which will be used in the following sections. For an elaborated treatment, one can refer to \cite{goubault2020directed}, \cite{grandis2003directed}.
\begin{definition} \label{d-space}
    A d-space is a topological space X with a distinguished set $dX \subseteq X^I$ of maps from the unit interval $I$ to $X$ (called d-paths) such that:
    \begin{itemize}
        \item $($Constant Paths$)$ All constant paths belong to $dX$. 
        \item $($Concatenation$)$ $dX$ is closed under concatenation.
        \item $($Composition$)$ $dX$ is closed under pre-composition with non-decreasing maps.
    \end{itemize}
\end{definition}
This will be referred to as a d-structure on $X$. We will denote a d-space $(X,dX)$ simply by $X$ when the context is clear. The set of d-paths $d X$ is a topological subspace of $ X^I$, equipped with the compact-open topology.
The subspace  \( \Gamma_X \subseteq X \times X \), is the collection of pairs \( (x, y) \in X \times X \) such that there exists a directed path in $dX$ from \( x \) to \( y \). We denote by $dX(x,x')$ the subspace of $dX$ consisting of d-paths from $x$ to $x'.$

A continuous map $f\colon X \to Y$ between two directed spaces $X$ and $Y$ is said to be a \emph{d-map} if it takes a d-path $\gamma$ in $X$ to a d-path $f \circ \gamma$ in $Y$. The induced map is denoted by $df\colon dX \to dY$. 

Denote by $\vv{I}$, the directed interval with the d-structure given by the collection of paths between $x$ and $y$ with $x \leq y$, for $x, y \in [0,1]$. Notice that for any d-space $X$ with a d-structure $dX$, we have $dX \subseteq X^{\vv{I}} \subseteq X^I$, where $X^{\vv{I}}$ is the space of d-paths on $X$. Note that, the largest d-structure on $X$ is given by the d-space $(X, X^I)$.

 Now we recall the definition of directed homotopy (see \cite{grandis2009directed}). A \emph{directed homotopy $($d-homotopy$)$} of d-maps $f,~g\colon X \rightarrow Y$  is given by a continuous map (known as a d-homotopy)
$$H\colon X \times \vv{I} \rightarrow Y$$
such that $H(x, 0)=f(x)$ and $H(x, 1)= g(x)$, for all $x \in X$. This is equivalent to a continuous map $H_t:=H(-, t)\colon X \rightarrow dY\subseteq Y^{\vv{I}}$ such that $H_0=f$ and $H_1=g$. When a pair of d-maps $f$ and $g$ are d-homotopic, we denote it by $f \simeq g$.

A pair of d-spaces are said to be d-homotopic if there exist d-maps $f\colon X \rightarrow Y$ and $f'\colon Y \to X$ such that $f' \circ f \simeq id_X$  and $f\circ f' \simeq id_Y$. We denote the d-homotopy inverses by $(f, f')$ for the d-spaces $X$ and $Y$.
\begin{remark}
Let  $(f, g)$ be d-homotopy inverses for d-spaces $X, Y$. Thus, the d-maps $f$ and $g$ induce maps
$$df\colon dX \rightarrow dY, \text{~and~} ~ dg\colon dY \rightarrow dX,$$
 such that  for any $x, x' \in X$, the restricted map is
 $$dg_{f(x), f(x')}\colon dY(f(x), f(x')) \rightarrow dX((g \circ f)(x), (g \circ f)(x')).$$
 Note that the map $dg_{f(x), f(x')}$ does not necessarily produce an element in $dX(x, x')$. 
 \end{remark}
 The above remark follows from the fact:
 Let $x, x' \in X$ and $\gamma \in dX((g \circ f)(x), (g \circ f)(x'))$. Since there is a d-homotopy $H\colon X \times \vv{I} \rightarrow X$ between $id_X$ and $g \circ f$, we get  the d-paths $\gamma_1(t):= H(x, t)$ and $\gamma_2(t):= H(x', t)$ of $dX(x, (g \circ f)(x))$ and $dX(x', (g \circ f)(x'))$ respectively. Also, we have a d-path  from  $(g \circ f)(x)$ to $(g \circ f)(x')$  given by $\gamma_3(t):=dg(H(x, t))$. 
 Since reversible path of a d-path is not a d-path, thus the usual way of concatenating paths $\gamma_1 * \gamma_3 * \gamma_2^{-1} $  to form a d-path from $x$ to $x'$ is not possible.
 Notice that $dg \circ df$ and $df \circ dg$ are not necessarily homotopic to $id_{dX}$ and $id_{dY}$ respectively.

As a result, we need a stronger notion of equivalence for d-spaces, namely \emph{basic dihomotopy equivalence} (see \cite{goubault2017directed, goubault2020directed}).

Recall that a topological space $X$ is said to be an Euclidean Neighborhood Retract (ENR) if it can be embedded into a Euclidean space $X \subset \mathbb{R}^N$ such that for an open neighborhood $X \subset U \subset \mathbb{R}^N$ there exists a retraction $r: U \rightarrow X$, $r|_X=id_X$. For example, any finite dimensional simplical complex is an ENR. Moreover, a subset $X \subseteq \mathbb{R}^n$ is an ENR if and only if it is locally compact and locally contractible, for $n \in \mathbb{N}$.

Suppose $X$ is an ENR.
Then, for a d-space $X$, the \emph{directed topological complexity}, denoted by $\tc(X)$, is defined as the smallest natural number $n$ such that there exists a cover $U_1,\dots, U_n$ of $\Gamma_X$ consisting of ENRs with each ENR admitting a continuous section of the dipath space map $\vv{\pi}\colon dX \to \Gamma_X.$
The directed topological complexity $\tc$ is a numerical invariant under such notion of equivalence (see \cite[Proposition 7]{goubault2020directed}).

We now recall the notion of a continuously graded map and then use it to introduce the concept of basic dihomotopy equivalence. 
\begin{definition}
    Let $f\colon X \to Y$ be a d-map. Let $y,y' \in Y$ and $W \subseteq X \times X$ be the inverse image of $(y,y')$ under the map $(f,f).$ Suppose we have continuous maps $$F^{y,y'}\colon dY(y,y') \times W \to dX$$ such that for all $(x,x') \in W$, we have $$F^{y,y'} (\gamma,x,x') \in dX(x,x').$$ In this case, we regard the map $F=(F^{y,y'})$ to be continuously graded. We also denote the grading of this map by $$F_{x,x'}: dY(f(x),f(x')) \to dX(x,x')$$ varying continuously over $(x,x') \in W$ in $dX^{dY(y,y')}$, with respect to the compact-open topology.
\end{definition}
\begin{definition}\label{defn: basic d-htopy equivalence}
    Let $X$ and $Y$ be a pair of d-spaces. A d-map $f\colon X \to Y$ is said to be a basic dihomotopy equivalence if the following conditions are satisfied:
    \begin{enumerate}
        \item $(f,g)$ is a d-homotopy inverses  between $X$ and $Y$.
        \item There exists a map $F\colon dY \to dX$ continuously graded by $$F_{x,x'}\colon dY(f(x),f(x')) \to dX(x,x') \text{~for~} x,x' \in \Gamma_X,$$ such that $(df_{x,x'}, F_{x,x'})$ is a homotopy equivalence between $dX(x,x')$ and $dY(f(x),f(x')).$
        \item There exists a map $G\colon dX \to dY$ continuously graded by $$G_{y,y'} : dX (g(y),g(y')) \to dY (y,y') \text{~for~} y,y'\in \Gamma_Y,$$ such that $(G_{y,y'},dg_{y,y'})$ is a homotopy equivalence between $dY(y,y')$ and $dX(g(y),g(y')).$
    \end{enumerate}    
\end{definition}
\begin{remark}
    In Definition \ref{defn: basic d-htopy equivalence}, if one of the spaces is a point, then other space is said to be dicontractible.
\end{remark}
Here, we recall the notion of directed fibrations (see \cite{grandis2002directed}, \cite{ioanpop}).
\begin{definition}
    A d-map $p\colon  E \to B$ is said to have d-homotopy lifting property with respect to a d-space $X$ if given d-maps $f\colon X \to E$ and $\varphi\colon X \times \vv{I} \to B$ and $\alpha \in \{0,1\}$ such that $\varphi \circ \partial^{\alpha} = p \circ f$, there is a directed lift $\varphi' \colon X \times \vv{I} \to E$ of $\varphi$ with respect to $p$ , $p\circ \varphi' =\varphi$ such that $\varphi' \circ \partial^{\alpha} = f.$ 
    \[\begin{tikzcd}
    X \arrow[r, "f"] \arrow[d, "\partial^{\alpha}"'] & E \arrow[d, "p"] \\
    X \times\vv{I} \arrow[r, "\varphi"] \arrow[ur, dashed, "\varphi'"]                 & B          
    \end{tikzcd}
    \] 
\end{definition}
\begin{definition}
    A d-map $p\colon  E \to B$ is called a d-fibration if it satisfies the d-homotopy lifting property with respect to every d-space. 
\end{definition}
\begin{remark}
    Given a d-fibration $p\colon  E \rightarrow  B$, we have $F_b \simeq F_{b'}$ for $b,b' \in B$. In other words, all fibres are d-homotopy equivalent, which we denote by $F$. 
\end{remark}
\begin{definition} \label{strong conn}
    A d-space $X$ is said to be strongly connected if $\Gamma_X= X \times X.$ In other words, there is a d-path between every pair of points of $X$.
\end{definition}
Now, we recall the notion of regular d-spaces, for which we have the product inequality: $\tc(X \times Y) \leq \tc(X)+ \tc(Y)-1$ for such spaces $X$ and $Y$ (see Proposition 1, \cite{goubault2020directed}).
\begin{definition}\label{d-regular}
    A d-space $X$ is called regular d-space if the directed end points space $\Gamma_X$ can be covered by $n$ ENRs as 
    $$\Gamma_X= A_1 \cup A_2 \cup \cdots \cup A_n, ~\text{where}~ \tc(X)=n$$
    with continuous sections over each $A_i$ of the dipath space map $\vv{\pi}:dX \rightarrow \Gamma_X$, 
    such that $A_i \cap A_j= \emptyset$ for $1 \leq i \neq j \leq n$ and the finite unions $A_1 \cup A_2 \cup \dots \cup A_r$ are closed for all $1 \leq r \leq n$. 
\end{definition}
\begin{remark}\label{rem: closure intersection regular d-space}
  We have the following property for a regular d-space with ENRs $\{A_i\}^n_{i=1}:$ 
    $$\bar{A_i} \cap A_j= \emptyset \text{~for $i < j$}.$$
\end{remark}

\section{Directed LS category and directed topological complexity}\label{sec: directed LS and directed TC}
In this section, we introduce the notion of directed LS category and study some of its properties. We also study its relationship with the directed topological complexity.
A motivation for defining this notion came from the Goubault's video seminar ``GEOTOP A: Directed topological complexity'' on CIMAT's youtube channel.

Recall that the Lusternik-Schnirelmann category of $X$, denoted by $\mathrm{cat}(X)$ is the least integer $n$ for which $X$ is covered by the $n$ contractible open sets in $X$.
In a classical setting it is known that $\mathrm{cat}(X)=\mathrm{secat}(p\colon P_0X\to X)$, where $p$ is the path fibration. There  are classical inequalities $\mathrm{cat}(X)\leq \TC(X)\leq \mathrm{cat}(X\times X)$ proved by Farber \cite[Theorem 5]{FarberTC}. We want to introduce the directed analogue of the Lusternik-Schnirelmann category so that this classical inequality also extends in the directed setting. 

A point $x_0\in X$ is called an \emph{initial point} if $(x_0,x)\in \Gamma_X$, for all $x\in X$.
Let $X$ be a d-space with initial point $x_0$ and consider the subspace of directed paths
$$d_0X:=\{\gamma\in dX \mid \gamma(0)=x_0\}.$$ Define an evaluation map $e_X\colon d_0X\to X$ by sending a d-path to its endpoint, i.e.,  $$e_X(\gamma):=\gamma(1).$$ We propose the following definition. 
\begin{definition} \label{d-cat}
The directed Lusternik-Schnirelmann (LS) category of a directed space $X$ is the smallest natural number $n$ (or infinity if it does not exist) for which $X$ can be covered by $n$ ENRs $U_1,\dots U_n$ such that each $U_i$ admits a continuous section of $e_X$.  
\end{definition}

\begin{remark} \label{motivation for directed cat}
 \normalfont{ In the directed settings, for a path-connected d-space $X$, existence of a continuous section of $\vv{\pi}\colon dX \rightarrow \Gamma_X$ is not equivalent to the dicontractibility of $X$ \cite[Theorem 1]{goubault2020directed}, unlike to the classical (undirected) case \cite[Theorem 1]{FarberTC}.
 We see this from the following \cite[Example 7]{goubault2020directed}: For $S^1$ with the  
  d-structure  given by any continuous path $\gamma\colon [0,1] \to S^1$ satisfies the following properties:
  (1) if $\gamma(0)=1$ then $\gamma$ is constant and 
  (2) if $\gamma(0) \neq 1$, then $| \gamma(t)+1 |$ is non-increasing. From this d-structure we obtain
 $$\Gamma_{S^1}=\{(e^{i\theta_1}, e^{i \theta_2})~|~ 0 < \theta_1 \leq \theta_2 < \pi\} \bigsqcup \{(e^{i\theta_1}, e^{i \theta_2})~|~ \pi < \theta_2 \leq \theta_1 < 2 \pi\} \bigsqcup \{(1,1), (-1,-1)\}.$$
 It follows that $\tc(S^1)=1$. But, $S^1$ is not a contractible d-space, thus not dicontractible.
 This implies that for an ENR $U\subseteq X$, unlike in the classical setting,  the existence of a continuous section of $e_X\colon d_0X\to X$ over $U$ may not imply the  dicontractibility of $U$ in $X$.
 Therefore, if we rewrite the above definition in terms of having minimal cover of dicontractible ENRs, we may not be able to establish the directed analogue of the inequality $\mathrm{cat}(X)\leq \TC(X)$. Nevertheless, we also note that the other implication, which is existence of continuous section of $\vv{\pi}$ whenever $U$ is dicontractible is always true. By \cite[Theorem 1]{goubault2020directed} there is a continuous section $s\colon \Gamma_U \to dU$.  }
\end{remark}
\[
\begin{tikzcd}
    dU\arrow[d, "\vv{\pi}|_U"]\\
    \Gamma_U \arrow[u, bend left, dashed,"s"]
\end{tikzcd}
~~~~~~~~\hspace{2 cm}
\begin{tikzcd}
    d_aU \arrow[r, "\phi"] & d_0X\\
    U\arrow[u,"s'"]\arrow[ur, sloped, "\sigma"]&
\end{tikzcd}\]
Here $s'(u)=s(a,u).$ Since $x_0$ is the initial point, there is a d-path joining $x_0$ to $a$, denoted by $\gamma_a$. Again, there is a d-path $s(a,u)$ joining $a$ to $u$. Let $\gamma_u \in d_a U$ and define $\phi(\gamma_u):=\gamma_a * \gamma_u$, obtaining a d-path joining $x_0$ to $u.$ Finally, define $\sigma:= \phi \circ s'.$

We now establish the basic dihomotopy invariance of directed LS category.

\begin{theorem}\label{thm: dihtopy equivalent then equal cat}
If $X$ and $Y$ are basic dihomotopy equivalent d-spaces with initial points, then $$\cat(X)=\cat(Y).$$  
\end{theorem}
\begin{proof}
Since $X$ and $Y$ are basic dihomotopy equivalent, there exists d-homotopy inverses $(f, g)$  for the d-spaces $X$ and $Y$. In addition, there exist continuously graded maps $F$ and $G$ such that $(df, F)$ and $(dg, G)$ are the homotopy inverses of $dX$ and $dY$.

Let $x_0$ and $y_0$ be the initial points for $X$ and $Y$, respectively. Define $d_0X = \{ \gamma \in dX ~|~ \gamma(0)=x_0\}$ and $d_0Y = \{ \gamma' \in dY ~|~ \gamma'(0)=y_0\}.$ Thus we have the induced maps $d_0f\colon d_0X \to d_0Y$ and $F_0 \colon d_0Y \to d_0X,$ which are d-homotopy inverses of each other.

Let $V \subseteq Y$ be an ENR with a continuous section $s'\colon V \rightarrow d_0 Y$ of $e_Y$. We want to show that the ENR $U:=f^{-1}(V)$ admits a continuous section $s: U \rightarrow d_0 X$ of $e_X$. To determine the section $s$ we use the following commutative diagram:
\[ \begin{tikzcd}
d_0 X \arrow{r}{d_0f} \arrow[swap]{d}{e_X} & d_0Y \arrow{d}{e_Y} \arrow{r}{F_0} & d_0 X \arrow{d}{e_X} \\
X \arrow{r}{f} & Y \arrow{r}{g} & X .
\end{tikzcd}
\]  
Our claim is that the required section $s$ is given by:
$$s(u) := \left( F_0 \circ s' \circ f \right) (u),$$ where $F_0(u):=F_{x_0,u}$ for all $u \in U$. For $v= f(u)$, we get $s'(v) \in d_0Y(y_0, v)$ where $y_0:=f(x_0)$, thus $F_0(s'(v)) \in d_0X(x_0, u)$ holds.
 Indeed we have, $e_X (s (u) )=e_X \left( F_0 (s' (v))\right) $. The grading $$F_{x_0,u}\colon dY(y_0,v) \to dX(x_0,u)$$ 
 and  $s'(v)\in dY(y_0,v),$ implies $F_0(s'(v))=F_{x_0,u} (s')(v) \in dX(x_0,u)$ and consequently under $e_X$ the d-path $F_0(s'(v))$ maps to $u.$ Thus, $\cat(X) \leq \cat(Y)$. Similarly, we can establish the reverse inequality.
\end{proof}
Now we state a comparison result for LS category and directed LS category. 
\begin{proposition}\label{prop: cat ineq & cat 1}
    Let $X$ be a d-space with an initial point. Then we have:
    \begin{enumerate}
        \item $\mathrm{cat}(X) \leq \cat(X).$
        \item If $X$ is dicontractible then $\cat(X)=1$. 
    \end{enumerate}
\end{proposition}
\begin{proof}
    (1) Suppose we have an ENR $U \subseteq X$, along with a continuous section $s\colon U \to d_0 X$ of $e_X$. Define a homotopy $H\colon X \times I \rightarrow X$ by
    $$H(x,t):=s(x)(1-t),$$ which gives the required nullhomotopy for $i_U: U \hookrightarrow X.$ Thus, the desired inequality follows using the above argument to a cover of $X$ with ENRs.

    (2) Since $X$ is dicontractible, it is basic dihomotopy equivalent to a point. Now applying \cref{thm: dihtopy equivalent then equal cat}, we obtain the desired result.
\end{proof}
\begin{remark}
We are unsure whether the converse to (2) of Proposition \ref{prop: cat ineq & cat 1} is true. To prove the converse, one needs to  reverse paths which is not allowed in the directed setting.
\end{remark}

We now establish the relationship between the directed LS category and the directed topological complexity.
\begin{proposition}\label{prop: rel cat-tc}
 Let $X$ be a d-space with an initial point $x_0$. Then 
\begin{enumerate}
    \item $\cat(X)\leq \tc(X) $.
    \item $\TC(X)\leq \cat(X\times X)$.
    \item If $X$ is regular, then $\TC(X)\leq 2\tc(X)-1$.
\end{enumerate}
\end{proposition}
\begin{proof}
$(1)$ Suppose $\tc(X)=n$. Then $\Gamma_X=\cup_{i=1}^{n} U_i$ with continuous sections $s_i\colon U_i\to dX$ of d-path space map.
Now define ENRs for $1\leq i\leq n$,  $$V_i:=\{v \in X ~|~ (x_0,v) \in U_i\} \text{~and~ } s'_i\colon V_i\to d_0X \text{~by~} s'_i(v):=s_i(x_0,v).$$
Then $s'_i$ is a continuous section of $e_{V_i}\colon  d_0V_i \to V_i.$ Indeed, for each $v \in V_i$, we have $(e_X \circ s'_i) (v) = e_X (s_i(x_0,v))= v$, implying $e_X \circ s'_i = id_{V_i}$. Note that the collection of ENRs $\{V_i\mid 1\leq i\leq n\}$ cover $X$.
This gives us the desired inequality $\cat(X) \leq \tc(X) $. 

$(2)$ Suppose $U\subseteq X\times X$ is an ENR with a section $s\colon U\to d_0(X\times X)$. 
Note that $X\times X$ has an initial point $(x_0,x_0)$ and thus $d_0(X\times X)=d_0X\times d_0X$. Then define the map $s'\colon U\to d_X$ by $$s'(x,y)= \left(pr_1\circ s(x,y)\right)^{-1}\ast \left(pr_2\circ s(x,y)\right),$$
 the concatenation of paths, the inverse path going from $x$ to $x_0$, with the d-path going from $x_0$ to $y$,
where $pr_1, pr_2\colon d_0X\times d_0X\to d_0X$ are the projections. Note that $s'(x,y)(0)=x$ and $s'(x,y)(1)=y$. Thus it defines a section of the free path space fibration $\pi\colon PX\to X\times X$.

$(3)$ Suppose $\tc(X)=n$.  Since $X$ is regular, we have a partition of $\Gamma_X=U_1\cup\dots \cup U_n$ into ENRs such that each $U_i$ admits a continuous section of $\pi\colon dX\to \Gamma_X$ and 
denote $\tilde{U}_i:=\{x\in X \mid (x_0,x)\in U_i\}$. Then note that the collection $\{\tilde{U}_i \mid 1\leq i\leq n\}$ covers $X$ and each $\tilde{U}_i$ admits a continuous section $s_i\colon \tilde{U}_i\to d_0X$ of $e_X$. Define $V_{i,j}=\tilde{U}_i\times \tilde{U}_j$ and $W_{k}=\cup_{i+j=k}V_{i,j}$ for $2\leq k\leq 2n$.  Then define $\sigma_{i,j}\colon V_{i,j}\to PX$ by $$\sigma(x,y)=(s_i(x))^{-1}\ast s_j(y),$$ where $(s_i(x))^{-1}$ is the inverse path. Note that $\sigma_{ij}$ defines a continuous local section of $\pi\colon PX\to X\times X$. Further note that $V_{ij}\cap V_{i'j'}=\emptyset$ for $(i,j)\neq (i',j')$. Thus we have a continuous section of $\pi$ over each $W_k$ (using Remark \ref{rem: closure intersection regular d-space}) for $2\leq k\leq 2n$. Moreover, the collection $\{W_k \mid 2\leq k\leq 2n\}$ covers $X\times X$. This implies $\TC(X)\leq 2n-1$. 
\end{proof}

\begin{remark}
 The proof of part $(3)$ of  Proposition \ref{prop: rel cat-tc} is inspired by Goubault's video seminar ``GEOTOP A: Directed topological complexity'' on CIMAT's youtube channel.  
\end{remark}

\begin{example} \label{d-sphere} Consider the directed $n+1$ dimensional cube $\vv{I}^{n+1}:= (\vv{I})^{n+1}$. Consider its boundary $\partial(\vv{I}^{n+1})$ as the directed sphere $S^n$, where the d-paths are those paths which are non-decreasing in each coordinate. Note that the origin in $\R^{n+1}$ is the initial point.
In \cite{borat2020directed}, Borat and Grant shown that $\tc(S^n)=2$ for all $n \geq 1$.
Thus, as an immediate application of Proposition \ref{prop: cat ineq & cat 1} and Proposition \ref{prop: rel cat-tc} we have $\cat(S^n)=2$.   
\end{example}
\begin{remark}
Observe the following equalities: $$\cat(S^{2n+1})=\tc(S^{2n+1})=\mathrm{cat}(S^{2n+1})=\TC(S^{2n+1})=2.$$
\end{remark}

\section{Directed Parametrized Topological Complexity} \label{Sec 4: directed parametrized topological complexity}
In this section, we introduce notion of directed parametrized topological complexity. 
Before we introduce this notion we begin by defining the parametrized topological complexity introduced by Cohen, Farber and Weinberger in \cite{C-F-W}.

For a fibration $p \colon E\to B$, consider the subspace $E^I_B$ of the free path space $E^I$ of $E$ defined by
\[
    E^I_B := \{\gamma\in E^I \mid \gamma(t) \in p^{-1}(b) ~\text{for some}~ b\in B ~\text{and for all}~ t\in[0,1] \}.
\]
Consider the pullback corresponding to the fibration $p \colon E\to B$ defined by 
\[
    E \times_B E = \{(e_1,e_2) \in E \times E \mid p(e_1)=p(e_2) \}.
\]
Then the \emph{parametrized endpoint map}
\begin{equation}
    \Pi \colon E^I_B \to E\times_B E, \quad \Pi(\gamma) := (\gamma(0),\gamma(1))
\end{equation}
is a fibration (see \cite[Appendix]{PTCcolfree}).
\begin{definition}
The  parametrized topological complexity of a fibration $p\colon E\to B$, denoted by $\mathrm{TC}[p\colon E\to B]$, is defined as the smallest natural number $n$ (infinity if it does not exists) such that $E\times_B E$ is covered by $n$ open sets $U_1,\dots, U_n,$ where each open set $U_i$ admits a continuous section of $\Pi$. 
\end{definition}

Now we define a parametrized analogue of directed topological complexity, for which we set up a few notation. Given a d-fibration $p\colon  E \to B$ 
with $E$ and $B$ as ENRs,
we define subspaces
$$dE_B:=\{ \gamma \in dE ~|~ p \circ \gamma~ \text{ is constant} \}$$
and $$\Gamma_{E,B}:= \{ (e_1,e_2) \in E \times_B E~|~ \exists ~ \gamma \in dE_B \text{~such that~} \gamma(0)=e_1, \gamma(1)=e_2 \}.$$
The \emph{directed parametrized endpoint map} is defined as 
\[
 \vv{\Pi}\colon dE_B\to \Gamma_{E,B}, \quad \vv{\Pi}(\gamma) := (\gamma(0),\gamma(1)).
\]

\begin{remark}
    For a d-fibration $p\colon  E \rightarrow B$, the subspace $dE_B \subseteq dE$ forms a d-structure on $E$. 
\end{remark}
\begin{definition}
    The directed parametrized topological complexity of a d-fibration $p\colon  E \to B$, denoted by $\tc[p\colon  E \to B]$, is defined as the smallest natural number $n$ (infinity if it does not exist) such that $\Gamma_{E,B}$ is covered by $n$ ENRs $U_1,\dots,U_n,$ where each ENR $U_i$ admits a continuous section of the map $\vv{\Pi}$.   
    
\end{definition}
\begin{remark}
When $B$ is a point, then $\Gamma_{E,B}=\Gamma_E$ and  $dE_B=dE$. Consequently, $$\tc[p\colon E \to B] =\tc(E).$$   
\end{remark}

\subsection{Properties} \label{properties}
In this part, we study some key features of directed parametrized topological complexity and its relationship with other standard numerical invariants.

We now analyze the behavior of directed parametrized topological complexity under the pullbacks of d-fibrations.
\begin{proposition}\label{prop:ineq: pullback and tc} 
Suppose $p\colon E\to B$ is a d-fibration and $B'\subseteq B$. If $p'\colon E'=p^{-1}(B')\to B'$ is the restricted d-fibration.  
Then we have $$\tc[p'\colon E' \to B'] \leq \tc[p\colon  E \to B] .$$
In particular, if $F=p^{-1}(b)$, then   $$\tc(F)\leq \tc[p\colon E\to B].$$
\end{proposition}
\begin{proof}
Observe that we have a commutative diagram:
 \[
    \begin{tikzcd}
    dE'_{B'}  \arrow[r, hook] \arrow[d, "\vv{\Pi'}"'] & dE_B \arrow[d, "\vv{\Pi}"] \\
    \Gamma_{E', B'} \arrow[r, hook]                                    & \Gamma_{E, B}
    \end{tikzcd}
    \]
where $\vv{\Pi'}$ is the restriction of $\vv{\Pi}$.
Suppose $U\subseteq \Gamma_{E, B}$ is an ENR with a continuous section $s\colon U\to dE_B$ of $\vv{\Pi}$. Define $V:=\Gamma_{E',B'}\cap U$. Note that $V$ is an ENR.
To obtain a continuous section of $\vv{\Pi'}$ on $V$, we will now show that the image of $s(V)$ lies inside $dE'_{B'}$. Suppose $(x,y)\in V$. Then $s(x,y)\in dE_B$. This implies 
$p( s(x,y)(t))=b$ for some $b\in B$ and for all $t\in \vv{I}$. 
In particular, $p( s(x,y)(t))=b=p( s(x,y)(0))=p(x)$. Since $x\in V\subseteq E'$, we have $b\in B'$. This gives us $s(x,y)(t)\in E'$ for all $t\in \vv{I}$. This shows that $s(V)\subseteq dE'_{B'}$. Then we can define $s':=s|_{V}\colon V\to dE'_{B'}$. This gives us the desired section of $\vv{\Pi'}$.

The inequality $\tc(F)\leq \tc[p\colon E\to B]$ follows by  setting $B'=\{b\}$.    
\end{proof}

We consider a parametrized version of  \cite[Theorem 1]{goubault2020directed} in the following.

\begin{theorem}\label{thm: dicontractibility implies tc 1}
Let $p\colon  E \to B$ be a d-fibration with fibre $F$ being either contractible d-space or it has an initial point. Then $F$ is dicontractible if and only if $$\tc[p\colon  E \to B]=1.$$
\end{theorem}
\begin{proof}
    Suppose $\tc[p\colon  E \to B]=1$.  Then from Proposition \ref{prop:ineq: pullback and tc}, we  have $\tc(F)=1$. 
     Now if $F$ has an initial point, then from Proposition \ref{prop: rel cat-tc} we have 
$\cat(F)=1$. This implies $F$ is contractible d-space. Then from  \cite[Theorem 1]{goubault2020directed}, we conclude that $F$ is dicontractible.

    Conversely, let $F$ be dicontractible. \cite[Theorem 1]{goubault2020directed} implies $\tc(F)=1$ and hence there exists a global continuous section $s$ of the d-map $p'\colon dF \to \Gamma_F.$ We have a  commutative diagram as follows: 
    \[\begin{tikzcd}
     dF \arrow[r, hook, "\alpha"] \arrow[ d, "p'"] & dE_B \arrow[d, "\vv{\Pi}"] \\
    \Gamma_F \arrow[u, bend left, "s"] \arrow[ r, hook ]                  & \Gamma_{E,B} \arrow[dashed, l, bend right, "\beta"]              
    \end{tikzcd}
    \]
    Choose $(e_1,e_2) \in \Gamma_{E,B}$. Then $p(e_1)=b= p(e_2).$ Denote $F_b : = p^{-1}(b).$ We have a d-homotopy $H_b \colon F_b \times \vv{I} \to F$  as all fibres are d-homotopic. Define: 
    $$\beta((e_1,e_2)) := (H_b(e_1,1), H_b(e_2,1))$$
    Since $(e_1,e_2) \in \Gamma_{F_b}$, we have $(H_b(e_1,1), H_b(e_2,1))\in \Gamma_F$, thus the map $\beta$ is well-defined. Finally, we have the composition $\alpha \circ s \circ \beta$, which is a global continuous section for $\vv{\Pi}$.
\end{proof}
\begin{proposition}\label{prop: trivial d-fib ; equality of tc}
For a trivial d-fibration $p\colon  E \to B$, we have $$\tc[p\colon  E \to B] =\tc(F).$$
\end{proposition}
\begin{proof}
Without loss of generality, consider $E=B \times F$, and suppose $p=pr_1\colon B \times F \to B$ and $pr_2\colon B \times F \to F$ are the first and second projections.
First we show that there are following canonical  d-homeomorphisms: $$\varphi\colon  \Gamma_{E,B} \xrightarrow{\cong} \Gamma_F \times B, ~ \psi: dE_B \xrightarrow{\cong} dF \times B.$$
To see these, define $\varphi\colon \Gamma_{E,B} \to \Gamma_F \times B$ by $(\tilde{e}_1, \tilde{e}_2) \mapsto ((e_1, e_2), b)$ where $pr_1(\tilde{e}_1)=b=pr_1(\tilde{e}_2)$ and $pr_2(\tilde{e}_1)=e_1$, $pr_2(\tilde{e}_2)=e_2$. Next define $\psi\colon dE_B \to dF \times B$ given by $\tilde{\gamma} \mapsto (\gamma, b)$, where $pr_1(\tilde{\gamma}(t))=b$ and $\gamma(t)=pr_2(\tilde{\gamma}(t)) \in F$ holds for all $t \in {\vv{I}}$.

Thus, it induces an isomorphism of the associated d-maps, given by  the following commutative  diagram
\[\begin{tikzcd}
 dE_B \arrow{d}{\vv{\Pi}} \arrow{r}{\psi} & dF \times B \arrow{d}{\vv{\pi} \times Id_B} \\%
\Gamma_{E, B} \arrow{r}{\varphi} & \Gamma_F \times B
& 
\end{tikzcd}\]
Note that Proposition \ref{prop:ineq: pullback and tc} already gives us one inequality $\tc(F) \leq \tc[p\colon  E \to B]$. In order to prove the other inequality, we begin with a partition of $\Gamma_F$ into ENRs $\{U_i\}_{i=1}^n$ and  sections $s_i$ over $U_i$ of the dipath space map $\vv{\pi}\colon dF \to \Gamma_F$ for $1\leq i\leq n$. Now consider the map $s_i \times Id_B: U_i \times B \to dF  \times B$. 
Define $V_i := \varphi^{-1} (U_i \times B)$ and $s_i':= \psi^{-1}\circ (s_i \times Id) \circ  \varphi$ for $1\leq i\leq n$. Note that $s_i'$ defines a continuous local section of $\vv{\Pi}$. Moreover, the collection $\{V_i\mid 1\leq i\leq n\}$ of ENRs cover $\Gamma_{E,B}$. This gives us the required inequality $\tc[p\colon  E \to B]=\tc[pr_1\colon B \times F \to B] \leq \tc(F)$. 
\end{proof}
\begin{corollary}\label{coro: tc of fiber and tc^bar for trivial d-fib}
    If a d-fibration $p\colon  E \rightarrow B$ is d-homotopic to the trivial d-fibration $pr_1:B \times F \rightarrow B$ for some path-connected space $F$, then $\tc[p\colon  E \to B] =\tc(F)$.
\end{corollary}

In what follows, we construct  d-fibrations $p\colon E \to B$ and $p'\colon E' \to B'$ satisfying $$\TC[p\colon E \to B] < \tc[p\colon E \to B]~~\text{~and~}~~\TC[p'\colon E' \to B'] > \tc[p'\colon E' \to B'].$$
\begin{example} \label{undirected TC < directed TC}
\normalfont{
    We endow $D^2$ with a particular d-structure (following \cite[Example 6]{goubault2020directed}) such that $$\TC(D^2)=1< \tc(D^2)=2.$$ The d-structure is given as follows : any directed path $\gamma \colon [0,1] \to D^2$ which starts at an internal point is constant. The directed path $\gamma\colon [0,1] \to D^2$, with $| \gamma(0)|=1$ are either of the two types : $\gamma(t)=e^{i\alpha(t)}$ or $\gamma(t)=e^{-i\alpha(t)}$, where $\alpha(t)$ is a non-decreasing function. 

    Consider the trivial fibration $p\colon B \times D^2 \to B$, with the specified d-structure on $D^2$, any d-structure on $B$ and $p$ is the projection map onto the first factor. By Corollary \ref{coro: tc of fiber and tc^bar for trivial d-fib}, we have $$\tc[p\colon B \times D^2 \to B]=\tc(D^2)=2.$$
    However, for the undirected case we have $$\TC[p\colon B \times D^2 \to B]=\TC(D^2)=1.$$}
\end{example}

\begin{example}\label{directed TC< undirected TC}
\normalfont{
    Consider $S^1$ with the following d-structure (following \cite[Example 7]{goubault2020directed}): any continuous path $\gamma\colon [0,1] \to S^1$ satisfies the following properties: (1) if $\gamma(0)=1$ then $\gamma$ is constant and (2) If $\gamma(0) \neq 1$, then $| \gamma(t)+1 |$ is non-increasing. For this d-structure we have $\tc(S^1)=1$. 

    Let $p'\colon B \times S^1 \to B$ be the trivial d-fibration projecting onto the first factor with $S^1$ being equipped with the specified d-structure and $B$ with any d-structure. Then we have $$\tc[p'\colon B \times S^1 \to B]=\tc(S^1)=1.$$ However in the undirected case we have  $$\TC[p'\colon B \times S^1 \to B]=\TC(S^1)=2.$$}
\end{example}

In \cite[Proposition 2]{goubault2020directed} it was shown that for strongly connected d-spaces, the topological complexity is bounded above by the directed topological complexity. 
Next we aim to establish the parametrized analogue of this result. First, we prove an essential lemma.
 
\begin{lemma}\label{lem: sc}
\begin{enumerate}
        \item Suppose $X, X'$ are basic dihomotopy equivalent d-spaces. Then $X$ is strongly connected if and only if $X'$ is strongly connected.
        \item Let $p\colon  E \to B$ be a d-fibration. Then $\Gamma_{E,B}=E \times_B E$, if the fibre $F$ of $p$ is strongly connected.
    \end{enumerate}
\end{lemma}
\begin{proof}
 (1)
 Suppose the basic dihomotopy equivalence between $X$ and $X'$ is given by the maps $f$ and $g$, which are d-homotopy inverses of each other. Start with a point $(e_1, e_2) \in X \times X$. Denote $e_1'=f(e_1)$ and $e_2'=f(e_2)$. Then $(e_1',e_2') \in X' \times X'$ and by the strongly connected hypothesis on $X'$, there exists a d-path $\gamma$ joining $e_1'$ and $e_2'$ in $X'$. Recall that, $(df_{e_1,e_2}, F_{e_1,e_2})$ is a homotopy equivalence between $dX(e_1,e_2)$ and $dX'(e'_1,e'_2),$ as in Definition \ref{defn: basic d-htopy equivalence}. Thus $F_{e_1,e_2}(\gamma)$ is the required d-path joining $e_1$ and $e_2$ in $X.$ The converse can be proved similarly.

(2) The inclusion $\Gamma_{E,B}\subseteq E\times_B E$ is obvious. To show the other inclusion, consider $(e_1,e_2) \in E \times_B E$. Suppose $e_1, e_2\in p^{-1}(b)=F_b$ where $p(e_1)=b=p(e_2)$. Since $F_b$ is strongly connected we have a directed path in $F_b$ whose end points are $e_1$ and $e_2$. 
Since all fibres are d-homotopic to each other, we get the desired equality.  
\end{proof}

The following result establishes the relationship between the parametrized topological complexity and its directed version.
\begin{proposition}\label{prop: dir and undir paramtrized TC inequality}
If the fibre $F$ of a d-fibration $p\colon E\to B$ is strongly connected, then we have
    \[
    \mathrm{TC}[p\colon E\to B]\leq \tc[p\colon E\to B].
    \]
\end{proposition}
\begin{proof}
Since $F$ is strongly connected, it follows from Lemma \ref{lem: sc} that $\Gamma_{E,B}=E\times_B E$. Therefore, any local continuous section of $\vv{\Pi}\colon dE_B \to \Gamma_{E,B}$ can be thought of as a local continuous section of $\Pi\colon E^I_B\to E\times_B E$. This completes the proof.
\end{proof}


\begin{definition}\label{def:  d secat regular}
    A d-fibration $p\colon  E \to B$ is called d-regular if $\Gamma_{E,B}$ can be covered by $n$ ENRs 
    $$\Gamma_{E,B}= A_1 \cup A_2 \cup \dots \cup A_n ~~\text{where}~~ \tc[p\colon E\to B]=n,$$
    with continuous local sections over each $A_i$ of $\vv{\Pi}:dE_B \rightarrow \Gamma_{E,B}$. Moreover, the ENRs satisfies  $A_i \cap A_j= \emptyset$ for $i \neq j$ and the finite unions $A_1 \cup A_2 \cup \dots \cup A_r$ are closed for all $1 \leq r \leq n$. 
\end{definition}
Note that for $B=\{\ast\}$, the notion of regular fibration recovers the notion of regular d-spaces.
\begin{remark}\label{rem: closure intersection regular d-map}
    The following property holds for the sets of a regular d-map with ENRs $\{A_i\}^n_{i=1}$ 
    $$\bar{A_i} \cap A_j= \emptyset \text{~for $i < j$. ~}$$
\end{remark}
We are now ready to prove the product inequality for directed parametrized topological complexity. 
\begin{proposition}\label{prop: product ineq dir parametrized tc}{}
    Let $p\colon  E \to B$ and $p'\colon E' \to B'$ be regular d-fibrations. Then we have: 
    $$\tc [ p \times p' \colon  E \times E' \to B \times B'] \leq \tc [p \colon E \to B] + \tc [ p' \colon  E' \to B']-1 .$$
\end{proposition}
\begin{proof}
Let $\tc [p \colon E \to B]=m$ and $\tc [p' \colon E' \to B']=n$.
 Thus, there are partitions by ENRs
     $$\Gamma_{E,B}= U_1 \cup U_2 \cup \dots \cup U_m ~~\text{and}~~ \Gamma_{E',B'} = V_1 \cup V_2 \cup \dots \cup V_n,$$ where we have continuous sections $s^1_i$ and $s^2_j$ on each of the ENRs $U_i$ and $V_j$ for the d-fibrations $p$ and $p'$ respectively. We consider the following collection of ENRs $G_r:= \bigcup\limits_{i+j=r} U_i \times V_j $ of $\Gamma_{E,B} \times \Gamma_{E',B'}$, for $1 \leq i \leq m$ and $1\leq j \leq n$. Clearly, the elements of $\{G_r\}_{r\geq 2}$ are pairwise disjoint  and form a cover of ENRs for $\Gamma_{E \times E', B \times B'} \cong \Gamma_{E,B} \times \Gamma_{E',B'}$.  
     Then $s^1_i\times s^2_j$ defines a continuous section of $\vv{\Pi}\times \vv{\Pi}'$ over $U_i\times V_j$. Then it follows that $U_i\times V_j$ is open in $G_r$, whenever $i+j=r$. Thus by Remark \ref{rem: closure intersection regular d-map}, the collection $$\bigsqcup\limits_{i+j=r} s^1_i\times s^2_j\colon G_r\to d(E\times E')$$ defines a continuous section of $\vv{\Pi}\times \vv{\Pi}'$ for $2\leq r\leq m+n$. This proves our desired inequality.
\end{proof}
We now present an inequality relating (undirected) parametrized topological complexity and directed LS category.
\begin{proposition}\label{prop: paramtrized tc leq directed LS cat} 
Let $p\colon  E \to B$ be a d-fibration and $E\times_B E$ has an initial point. Then 
 $$\mathrm{TC}[p\colon  E \to B] \leq \cat(E \times_B  E).$$
\end{proposition}
\begin{proof}
Note that, for the underlying fibration, we have $\TC[p\colon E\to B]\leq \mathrm{cat}(E\times_B E)$. Then the desired inequality follows from Proposition \ref{prop: cat ineq & cat 1}.   
\end{proof}

\subsection{Invariance}\label{invariance}
In this subsection, we establish the fibrewise basic dihomotopy invariance of the directed parametrized topological complexity. 
To introduce the fibrewise basic dihomotopy equivalence of d-fibrations, we      
need the following definitions.
\begin{definition}
Let $p\colon E\to B$ and $p'\colon E'\to B$ be d-fibrations. 
A fibrewise map from $p\colon E\to B$ to $p'\colon E'\to B$ is a map $f\colon E\to E'$ such that $p'\circ f=p$.    
\end{definition}

\begin{definition}
A fibrewise d-homotopy $F\colon E\times \vv{I}\to E'$ is a map such that $q(F(-,t))=p$ for all $t\in \vv{I}$. Thus, $F$ is a d-homotopy between fibrewise maps $F(-,0)$ and $F(-,1)$.
\end{definition}
\begin{definition}\label{def: d-fib hteq} 
The d-fibrations $p\colon E\to B$ and $p'\colon E'\to B$ are said to be fibrewise basic dihomotopy equivalent if the following conditions hold:
\begin{enumerate}
    \item  There exist fibrewise maps $f\colon E\to E'$ and $g\colon E'\to E$ such that there are fibrewise d-homotopies  from $f\circ g$ to $Id_{E'}$ and from $g\circ f$ to $Id_{E}$.
    \item  There exist a continuously graded map $F: dE' \rightarrow dE$ satisfying the commutative diagrams
\[\begin{tikzcd}[every arrow/.append style={shift left}]
dE \arrow[dr, "dp"']  \arrow[rr, "df"  shift left=1.5ex] 
& & dE' \arrow[dl, "dp'"]  \arrow[ll, "F"  shift right=1.5ex] \\
& dB 
\end{tikzcd}\]
 such that for $(e_1, e_2) \in \Gamma_E$, we have gradings $F_{e_1, e_2}: dE'(f(e_1), f(e_2)) \rightarrow dE(e_1, e_2)$ and a homotopy equivalence $(df_{e_1, e_2}, F_{e_1, e_2})$  between $dE'(f(e_1), f(e_2))$ and $dE(e_1, e_2)$.
 \item There exist a continuously graded map $G: dE \rightarrow dE'$ satisfying the commutative diagrams
\[\begin{tikzcd}[every arrow/.append style={shift left}]
dE' \arrow[dr, "dp'"']  \arrow[rr, "dg"  shift left=1.5ex] 
& & dE \arrow[dl, "dp"]  \arrow[ll, "G"  shift right=1.5ex] \\
& dB 
\end{tikzcd}\]
such that for $(e'_1, e'_2) \in \Gamma_{E'}$, we have gradings $G_{e'_1, e'_2}: dE(g(e'_1), g(e'_2)) \rightarrow dE'(e'_1, e'_2)$ and a homotopy equivalence $(dg_{e'_1, e'_2}, G_{e'_1, e'_2})$  between $dE(g(e'_1), g(e'_2))$ and $dE'(e'_1, e'_2)$.
\end{enumerate}
\end{definition}

\begin{example}
Suppose $p\colon  E\to B$, $p'\colon E'\to B$ are d-fibrations and $f$ is a d-homeomorphism satisfying the following commutative diagram    
\[\begin{tikzcd}[every arrow/.append style={shift left}]
E \arrow[dr, "p"']  \arrow[rr, "f"] 
& & E' \arrow[dl, "p'"]  \\
& B.
\end{tikzcd}\]
Then $p$ and $p'$ are fibrewise basic dihomotopy equivalent. Note that the data $f$, $g=f^{-1}$, $dg=F$, and $df=G$ establishes the fibrewise basic dihomotopy equivalence of $p$ and $p'$.
\end{example}
\begin{remark}
  In Definition \ref{def: d-fib hteq}, substituting the base space as a point, we recover Definition \ref{defn: basic d-htopy equivalence}.  
\end{remark}

With the previous definitions, we establish the fibrewise basic dihomotopy invariance of directed parametrized topological complexity.
\begin{theorem}\label{thm: basic dihomotopy eqv fib have same TC}
If d-fibrations $p\colon E\to B$ and $p'\colon E'\to B$ are fibrewise basic dihomotopy equivalent, then $$\tc[p\colon E\to B]=\tc[p'\colon E'\to B].$$   
\end{theorem}
\begin{proof}
Since $p\colon E\to B$ and $p'\colon E'\to B$ are fibrewise basic dihomotopy equivalent, we have d-maps $f\colon E \to E'$ and $g\colon E' \to E$  which are d-homotopy inverses of each other satisfying the following commutative diagram:  
\[\begin{tikzcd}[every arrow/.append style={shift left}]
E \arrow[dr, "p"']  \arrow[rr, "f" shift left=1.5ex] 
& & E' \arrow[dl, "p'"]  \arrow[ll, "g" shift right=1.5ex] \\
& B.
\end{tikzcd}\]
Moreover, there is a d-map $G\colon dE\to dE'$ satisfying part $(3)$ of Definition \ref{def: d-fib hteq}. Therefore, $G$ restricts to a map from $dE_B$ onto $dE'_B$, which we again denote by $G$. Thus, we have the following  commutative diagram

\[ \begin{tikzcd}
dE'_B \arrow{r}{dg} \arrow[swap]{d}{\vv{\Pi'}} & dE_B \arrow{d}{\vv{\Pi}} \arrow{r}{G} & dE'_B \arrow{d}{\vv{\Pi'}} \\%
\Gamma_{E', B} \arrow{r}{\tilde{g}}& \Gamma_{E, B}  \arrow{r}{\tilde{f}} & \Gamma_{E', B}
\end{tikzcd}
\] 
where the d-map $\tilde{g}$ is defined as $\tilde{g}:=g\times g|_{\Gamma_{E',B}}$ and $\tilde{f}$ defined similarly. Note that these maps are well defined as they are fibrewise maps. Let's verify this for completeness.
Let $(e_1',e_2') \in \Gamma_{E',B}.$ Then we have $p'(e_1')=b=p'(e_2')$ and consequently $p(g(e_1'))=b=p(g(e_2')),$ from which it follows that $(g(e_1'), g(e_2')) \in \Gamma_{E,B}.$ We also check that the map $G\colon dE \to dE'$ on $dE_B$ restricts to $dE_B \to dE'_B.$  For $\gamma \in dE_B$, we have $(dp' \circ G) (\gamma)=dp(\gamma)$. Thus, $p'(G(\gamma)(t))= p(\gamma(t))=b$ for all $t \in \vv{I}$, and hence $G(\gamma) \in dE'_B$. 

Let $U\subseteq \Gamma_{E, B}$ be an ENR such that there is a continuous section $s$ of $\vv{\Pi}$ on $U$.
Define $V:=\tilde{g}^{-1}(U)$
and $s'\colon V\to dE'_B$ by $s'(v):= G_v \circ s \circ \tilde{g}|_{V}(v)$ for $v\in V$. 
For $(e'_1, e'_2) \in \Gamma_{E', B}$, the path  $\gamma:=s((g(e'_1), g(e'_2))) \in dE_B$, i.e. $\gamma \in dE_B(g(e'_1), g(e'_2))$.
Since, $G(\gamma) \in dE'_B(e'_1, e'_2)$, we have $\vv{\Pi'}(G(\gamma))=(e'_1, e'_2)$. Note that $s'$ is continuous as $s$ is continuous, $\tilde{g}$ is continuous and $G$ is continuous and continuously graded.
Applying the above argument to cover with ENRs, we obtain the following inequality:
$$\tc[p'\colon E'\to B] \leq \tc[p\colon E\to B].$$

Now we can apply the same argument to the following diagram to obtain the reverse inequality.
\[ \begin{tikzcd}
dE_B \arrow{r}{df} \arrow[swap]{d}{\vv{\Pi}} & dE'_B \arrow{d}{\vv{\Pi'}} \arrow{r}{F} & dE_B \arrow{d}{\vv{\Pi}} \\%
\Gamma_{E, B} \arrow{r}{\tilde{f}}& \Gamma_{E', B}  \arrow{r}{\tilde{g}} & \Gamma_{E, B}.
\end{tikzcd}
\] 
This completes the proof.
\end{proof}

\section{Some computations} \label{Sec 5: computations}
In this section we consider a particular kind of non-trivial but very natural directed structure on a fibration. Using that we  compute the directed parametrized topological complexity for the Hopf fibration and the Fadell-Neuwirth fibration.
\begin{definition} \label{special d-str on E}
   Let $p\colon E \rightarrow B$ be a fibration with a directed structure $dB$ on $B$.
Define the new directed structure $dE$ on $E$ as follows:
   $$\gamma \in dE ~\text{if and only if}~ p \circ \gamma \in dB.$$
\end{definition}

Note that for a path $\gamma$ in $F_b:=p^{-1}(b)$, i.e. $\gamma \in E^I_B$, $p \circ \gamma$ is a constant path at $b$.
\begin{remark}\label{d-struc on fibre}
With such a d-structure on the total space $E$ of a fibration $p\colon E \to B$, on each fibre $F_b$ we get an induced d-structure from $E$, for all $b \in B$. Without loss of generality, we denote a fibre by $F$, together with the d-structure.
For any path $\gamma$ in $F$, $p\circ \gamma$  is a constant path. Therefore, a path $\gamma$ automatically becomes a directed path with respect to the induced d-structure. As a result, this d-structure on $F$ is the largest one in the sense of Grandis \cite{grandis2003directed}. That is, this d-structure is nothing but the free path space of $F$, denoted by $F^I$. Consequently, $F$ is strongly connected and it follows that  $\tc(F)=\mathrm{TC}(F)$.
Also, observe that in this situation  $p\colon E \rightarrow B$ becomes a d-fibration.
\end{remark}

\begin{theorem} \label{directed-usual equality}
Suppose $E$ and $B$ are and $p\colon E\to B$ is a d-fibration as constructed above. Then 
    \[\tc[p\colon E\to B] = \mathrm{TC}[p\colon E\to B].\]   
\end{theorem}
\begin{proof}
Consider the usual parametrized endpoint map $\Pi \colon E^I_B \to E\times_B E$ associated to the fibration $p$. Suppose $\mathrm{TC}[p\colon E\to B]=n$. Since both $E$ and $B$ are ENRs, then there exists a cover of $E\times_B E$ consisting ENRs $\{U_1,\dots, U_n\}$ with sections $s_i$ of the fibration $\Pi$ on $U_i$. Our aim is to show that these sections are directed.
Note that $s_i(e_1, e_2) \in F_b$ holds, where $p(e_1)=b$ and $(e_1, e_2) \in U_i$ for  $i=1, \dots,n$. Since the induced d-structure on $F$ from the associated d-structure on $E$ for $p$ is given by the free path space $F^I$, $s_i(e_1,e_2)$ can be considered a directed path between $e_1$ and $e_2$. More specifically, $s_i$ is a directed map. Also note that, the fibre of $p\colon E\to B$ are strongly connected. Therefore, from part-2 of Lemma \ref{lem: sc}, we have $\Gamma_{E,B}=E\times_B E$. Then the sections $s_i$ on $U_i$'s can be considered as directed sections of  $\vv{\Pi}\colon dE_B\to \Gamma_{E,B}$ for $1\leq i\leq n$. This implies $\tc[p\colon E\to B] \leq \TC[p\colon E\to B]$.
The other inequality follows from Proposition \ref{prop: dir and undir paramtrized TC inequality} as $F$ is strongly connected.
\end{proof}

\begin{corollary} \label{TC(F)-TC[p]} Let $p\colon E\to B$ be a d-fibration as mentioned above. Thus
    by Proposition \ref{prop:ineq: pullback and tc} and by Remark \ref{d-struc on fibre}, we get  the following inequality
    $$\mathrm{TC}(F)\leq \tc[p\colon E\to B].$$
\end{corollary}
Using the above results we evaluate the directed parametrized topological complexity of the directed fibrations associated with the Hopf fibrations and the Fadell-Neuwirth fibrations.

\subsection{Directed Hopf fibrations} Consider the Hopf fibration
$$S^1 \hookrightarrow S^{3} \xrightarrow{p} S^2.$$ This fibration also has a principal $S^1$-bundle structure. Using the standard d-structure on $S^{2}$ as mentioned in Example \ref{d-sphere}, we can give an associated d-structure on $S^{3}$ such that the above fibration becomes a directed one. Since in the usual case, we have $\mathrm{TC}[p\colon  S^{3} \rightarrow S^{2}]=2$ (see \cite[Example 4.4]{C-F-W} or \cite[Proposition 4.3]{C-F-W}), together with \cref{directed-usual equality} we get
$$\tc[p\colon S^{3} \rightarrow S^{2}] = 2.$$

For the general Hopf fibration $S^1 \hookrightarrow S^{2n+1} \xrightarrow{p} \mathbb{C}P^{n}$ (which are  principal $S^1$-bundles), we can again use \cref{directed-usual equality} and  \cite[Proposition 4.3]{C-F-W} to conclude that 
$$
\tc[p\colon S^{3} \rightarrow S^{2}]=\TC[p\colon S^{2n+1} \xrightarrow{p} \mathbb{C}P^{n}]=\TC(S^1)=2.
$$
\subsection{Directed Fadell-Neuwirth fibrations} The problem of obstacle-avoiding, collision-free motion of multiple robots in the presence of multiple obstacles whose positions are not known in advance, was studied in \cite{C-F-W}, \cite{PTCcolfree}. To deal with this problem, the notion of \emph{Fadell-Neuwirth fibration} is an appropriate model.

Note that for a topological space $Y$, the configuration space of $n$ distinct ordered points in Y is denoted by
\[
F(Y, n) = \left\{ (x_1, \ldots, x_n) \in Y^n \,\middle|\, x_i \neq x_j \text{ for } i \neq j \right\}.
\]
We consider the \emph{Fadell--Neuwirth fibration} 
\begin{equation}\label{Fadell-Neuwirth}
    p : F(\mathbb{R}^k, n + m) \rightarrow F(\mathbb{R}^k, m),
\end{equation}

given by
\[
(z_1, \dots, z_m, z_{m+1}, \dots, z_{m+n}) \mapsto (z_1, \dots, z_m),
\]
and aim to compute the directed counterpart of its \emph{parametrized topological complexity}.

The parametrized topological complexity of these fibrations were computed in \cite{C-F-W} and \cite{PTCcolfree}. 
In Example \ref{d-sphere}, we mentioned the standard d-structure on the $n$-dimensional cube $\vv{I}^n$, which naturally extends to $\mathbb{R}^n$. This extension induces a canonical d-structure on the configuration space $F(\mathbb{R}^k, m)$, the base of the Fadell–Neuwirth fibration (\ref{Fadell-Neuwirth}). By equipping the total space $F(\mathbb{R}^k, n + m)$ with the associated d-structure as in Definition \ref{special d-str on E}, we obtain a directed version of the Fadell–Neuwirth fibration. Applying Theorem \ref{directed-usual equality} to this d-fibration yields:

\begin{theorem}
\label{thm: dptc-Fadell_Neuwirth}
Suppose $m\geq 2$, $n\geq 1$. 
Then
\[
    \tc[p \colon F(\R^k,m+n)\to F(\R^k, m)] =
\begin{cases}
    2n+m, & \text{if $k$ is odd},\\
    2n+m-1 & \text{if $k$ is even}.
\end{cases}
\]
\end{theorem}






\vspace{.5 cm}
\noindent\textbf{Acknowledgment.}
The authors would like to thank Prof. Goubault for his valuable suggestions. Moreover, Navnath Daundkar acknowledges the support of the DST–INSPIRE Faculty Fellowship (Faculty Registration No. IFA24-MA218), Department of Science and Technology, Government of India.
Abhishek Sarkar acknowledges the support of IISER Pune for the Institute Post-Doctoral fellowship IISER-P/Jng./20235445.

\bibliographystyle{plain} 
\bibliography{references}

\end{document}